\documentclass{tac}
\usepackage{epsfig,epic,eepic,latexsym, amssymb}
\copyrightyear{2006}

\title{Frobenius Objects in Cartesian Bicategories}
\author{R.F.C Walters and R.J. Wood}
\thanks{The authors gratefully acknowledge financial support from 
the Italian CNR and the Canadian NSERC.
Diagrams typeset using M. Barr's diagram package, diagxy.tex.}
\address{\\Dipartimento di Scienze delle Cultura\\
Politiche e dell'Informazione\\
Universit\`a dell Insubria, Italy\\[3pt] 
and\\[3pt]
Department of Mathematics and Statistics\\
Dalhousie University\\
Halifax, NS, B3H 3J5 Canada}

\eaddress{robert.walters@uninsubria.it, rjwood@dal.ca}

\newtheorem{axiom}{Axiom}

\let\thm\theorem
\let\prp\proposition

\let\lem\lemma
\let\dfn\definition
\let\eth\endtheorem
\let\prf\proof
\let\frp\endproof

\let\axm\axiom

\input diagxy
\xyoption{curve}

\newcommand{\s}{\scalefactor{0.5}}

\newcommand{\f}{{\kern -.25em}:{\kern -.25em}}

\newcommand{\ra}{{\s\to}}
\newcommand{\la}{\,{\s\toleft}\,}

\newcommand{\rass}{\s{|{\kern -.38em}\to}}

\newcommand{\iso}{\cong}
\newcommand{\laj}{\dashv}

\newcommand{\x}{\times}
\newcommand{\ox}{\otimes}
\newcommand{\Hat}{\widehat}
\newcommand{\Tilde}{\widetilde}

\newcommand{\op}{{^\mathrm{op}}}
\newcommand{\rev}{{^\mathrm{rev}}}
\newcommand{\oprev}{{^\mathrm{oprev}}}
\newcommand{\inv}{^{-1}}

\newcommand{\CAT}{\mathbf{CAT}}

\newcommand{\prof}{\mathbf{prof}}

\newcommand{\grp}{\mathrm{Frob}}
\newcommand{\map}{\mathrm{Map}}

\newcommand{\rel}{\mathrm{Rel}}
\newcommand{\spn}{\mathrm{Span}}

\newcommand{\E}{{\cal E}}

\newcommand{\bB}{\mathbf{B}}

\newcommand{\bG}{\mathbf{G}}
\newcommand{\bM}{\mathbf{M}}

\begin{document}

\maketitle

\begin{abstract}
Maps (left adjoint arrows) between Frobenius objects
in a cartesian bicategory $\bB$ are precisely comonoid homomorphisms
and, for $A$ Frobenius and any $T$ in $\bB$, $\map(\bB)(T,A)$ is
a groupoid.
\end{abstract}

\section{Introduction}\label{intro}
The notion of locally ordered cartesian bicategory was introduced by
Carboni and Walters \cite{caw} for the axiomatization of the bicategory 
of relations of a regular category. The notion has since been extended 
by Carboni, Kelly, Walters, and Wood \cite{ckww} to the case of 
a general bicategory, to include examples such as bicategories of 
spans, cospans, and profunctors.

A crucial further axiom introduced by Carboni and Walters in that paper
was the so-called discreteness axiom, now known as the Frobenius axiom,
since it was recognized to be equivalent to Lawvere's equational 
version \cite{law} of Frobenius algebra. With this axiom one can 
define the notion of Frobenius object in a
monoidal category, the Frobenius axiom being an equation satisfied by
monoid and comonoid structures on the object.

The Frobenius axiom has found a large variety of uses. For example, 
the 2-dimensional cobordism category has been shown to be the 
symmetric monoidal category with a generic commutative Frobenius 
object. (For a presentation of this result see J. Kock \cite{ko}.) 
Related results are the characterization of the symmetric monoidal 
category of cospans of finite sets in \cite{lack} and the 
characterization of the symmetric monoidal category of cospans of 
finite graphs in \cite{rsw}.
Another example is that, in the algebra of quantum measurement 
\cite{cp}, classical data types are Frobenius objects. In \cite{gah} 
the Frobenius equation is a crucial equation in an algebraic 
presentation of double pushout graph rewriting, and in \cite{ksw} 
the equation is one of the main equations in a compositional theory 
of automata. The 2-dimensional version of Frobenius algebra has also 
been introduced in the characterization of a certain monoidal 
2-category in \cite{msw}.

There is a rather obvious way of extending the notion of Frobenius 
object to the context of a monoidal bicategory: instead of requiring 
equations between operations, certain canonical 2-cells are required 
to be invertible. This paper develops properties of such 2-dimensional 
Frobenius objects, for the canonical monoid and comonoid structure 
on each object which is part of the cartesian bicategory structure. 
The two principal results are (i) that maps (left adjoint arrows)
between Frobenius objects are the same as comonoid homomorphisms, 
and (ii) that if $A$ is a Frobenius object then, for any object $T$
in the cartesian bicategory $\bB$, $\map(\bB)(T,A)$ is a groupoid.
This second result was noticed for the special case of Profunctors 
at the time of the Carboni-Walters paper by Carboni and Wood, 
independently, but has never been published. We develop in this paper 
techniques in a general cartesian bicategory which enable us to lift 
the profunctor proof.

The results of this paper will be used in a following paper \cite{waw} 
characterizing bicategories of spans.

We thank Bob Par\'e for his helpful comments when these results were 
presented at the ATCAT seminar in Halifax, Canada.

\section{Preliminaries}\label{prelim}
\subsection{}\label{pl1}
We recall from \cite{ckww} that a bicategory $\bB$ is 
{\em cartesian\/} if the subbicategory of maps (by which we mean
left adjoint arrows), $\bM=\map\bB$, 
has finite products ($-\x-$,~$1$) with projections denoted
$p\f X\la X\x Y\ra Y\f r$; each hom-category
$\bB(X,A)$ has finite products ($-\wedge-$, $\top$) 
with projections denoted $\pi\f R\la R\wedge S\ra S\f \rho$; 
and an evident derived tensor product on $\bB$, ($-\ox-$, $I$) 
extending the product structure of $\bM$, is functorial. 
It was shown that the derived tensor product of a 
cartesian bicategory underlies a symmetric monoidal bicategory 
structure. Throughout this paper, $\bB$ is assumed to be a 
cartesian bicategory and, as in \cite{ckww}, we assume, for
ease of notation, that $\bB$ is normal, meaning that the identity
compositional constraints of $\bB$ are identity 2-cells.

\subsection{}\label{pl2}
If $f$ is a map of $\bB$, an arrow of $\bM$, we will write 
$\eta_f,\epsilon_f\f f\laj f^*$ for a chosen adjunction in $\bB$ 
that makes it so. We occasionally refer to an $f^*$ as a {\em pam}.
As in \cite{ckww}, we write 
$$\bfig
\Atriangle/->`->`/[\bG`\bM`\bM;\partial_0`\partial_1`]
\efig$$
for the Grothendieck span corresponding to
$$\bM\op\x\bM\to^{i\op\x i}\bB\op\x\bB\to^{\bB(-,-)}\CAT$$
where $i\f\bM\ra\bB$ is the inclusion. A typical arrow of $\bG$,
$(f,\alpha,u)\f(X,R,A)\ra(Y,S,B)$ can be depicted by a square in $\bB$
\begin{equation}\label{Gsquare}
\bfig
\square(0,0)[X`Y`A`B;f`R`S`u]
\morphism(125,250)|a|<250,0>[`;\alpha]
\efig
\end{equation}
in which $f$ and $u$ are maps, and such arrows are 
composed by pasting. A 2-cell 
$(\phi,\psi)\f(f,\alpha,u)\ra(g,\beta,v)$ in $\bG$ is a pair of
2-cells $\phi\f f\ra g$, $\psi\f u\ra v$ in $\bM$ which satisfy
the obvious equation.

\subsection{}\label{pl3}
In part of this and subsequent work it will be useful to revisit
certain of the arrows of $\bG$ from another point of view. Consider
$$\bfig
\square(0,0)[T`T`X`Y;1_T`x`y`R]
\morphism(250,375)|r|<0,-250>[`;\rho]
\efig$$
On the one hand it is just an arrow from $1_T$ to $R$ in $\bG$
but each of the three reformulations of $\rho$ that result from
taking mates have their uses.
$$\bfig
\square(0,0)/->`->`<-`->/[T`T`X`Y;1_T`x`y^*`R]
\morphism(250,375)|r|<0,-250>[`;\Hat\rho]
\square(1000,0)/->`<-`<-`->/[T`T`X`Y;1_T`x^*`y^*`R]
\morphism(1250,375)|r|<0,-250>[`;\rho^*]
\square(2000,0)/->`<-`->`->/[T`T`X`Y;1_T`x^*`y`R]
\morphism(2250,375)|r|<0,-250>[`;\Tilde\rho]
\efig$$
In the first of these, $\Hat\rho\f 1_T\ra y^*Rx$, it is sometimes
convenient to write $R(y,x)=y^*Rx$ and regard $\Hat\rho$ as a
$1_T$-element of $R(y,x)$. In the special case where $R$ is
$1_X\f X\ra X$ we write $X(y,x)=y^*x$ (invoking normality of $\bB$).
(This hom-notation is similar to that employed first in \cite{saw}. 
It was adapted for this compositional context in \cite{wd1}.)
The second we will use without further comment except to say that,
for $R=1_X$, $\rho^*$ is the usual way of making the process
of taking right ajoints functorial. The third will appear in our
discussion of tabulations in the forthcoming \cite{waw}. Note that
the $R(y,x)$ notation extends to 2-cells so that, for $\eta\f y'\ra y$
and $\xi\f x\ra x'$, we have $R(\eta,\xi)\f R(y,x)\ra R(y',x')$.

The chief purpose of the notation $R(y,x)$ is to guide intuition
so that constructions in such cartesian bicategories as that 
of categories, profunctors, and equivariant 2-cells (which we call
$\prof$) can be usefully generalized. Observe that if 
$\tau\f R\ra S$ is a 2-cell in $\bB$ and $\xi\f x\ra x'$ then
we have automatically such identities as 
$\tau(y,x').R(y,\xi)=S(y,\xi).\tau(y,x)$, both providing the
horizontal composite $\tau\xi$ whiskered with $y^*$ as below.
$$\bfig
\morphism(0,0)|a|/{@{->}@/^2em/}/<500,0>[T`X;x]
\morphism(0,0)|b|/{@{->}@/_2em/}/<500,0>[T`X;x']
\morphism(500,0)|a|/{@{->}@/^2em/}/<500,0>[X`Y;R]
\morphism(500,0)|b|/{@{->}@/_2em/}/<500,0>[X`Y;S]
\morphism(250,100)|r|<0,-200>[`;\xi]
\morphism(750,100)|r|<0,-200>[`;\tau]
\morphism(1000,0)|a|<500,0>[Y`T;y^*]
\efig$$
For the most part, we will use such calculations with little
comment.

If 
$$\bfig
\square(0,0)/->`->`<-`->/[T`T`X`Y;1_T`x`y^*`R]
\morphism(250,375)|r|<0,-250>[`;\Hat\rho]
\place(750,500)[\mbox{and}]
\square(1000,0)/->`->`<-`->/[T`T`Y`Z;1_T`y`z^*`S]
\morphism(1250,375)|r|<0,-250>[`;\Hat\sigma]
\efig$$
are $1_T$-elements of $R(y,x)$ and $S(z,y)$ respectively then it is
easy to see that $\Hat{\rho\Box\sigma}$, where $\rho\Box\sigma$ is the
paste composite of $\rho$ and $\sigma$, is a $1_T$-element of 
$(SR)(z,x)$. The $1_T$-element $\Hat{\rho\Box\sigma}$ can be given
in several ways. We will have occasion to give it via the pasting 
composite
$$\bfig
\square(0,0)/->`->``/[T`T`X`;1_T`x``]
\morphism(0,0)|b|<750,0>[X`Y;R]
\morphism(300,375)|l|<0,-250>[`;\rho]
\Vtriangle(500,0)/->`->`<-/<250,500>[T`T`Y;1_T`y`y^*]
\morphism(750,400)|m|<0,-250>[`;\eta_y]
\square(1000,0)/->``<-`/[T`T``Z;1_T``z^*`]
\morphism(750,0)|b|<750,0>[Y`Z;S]
\morphism(1200,375)|r|<0,-250>[`;\sigma^*]
\efig$$

We note that a paste composite such as $\rho\Box\sigma$ as below
$$\bfig
\square(0,0)|almb|[T`T`X`Y;1_T`x`y`R]
\morphism(250,375)|l|<0,-250>[`;\rho]
\square(500,0)|amrb|[T`T`Y`Z;1_T`y`z`S]
\morphism(750,375)|r|<0,-250>[`;\sigma]
\efig$$
may result from several different $y\f T\ra Y$. For example, in
$$\bfig
\square(0,0)|almb|/->`->`{@{->}@/_2em/}`->/[T`T`X`Y;1_T`x`y`R]
\morphism(175,375)|l|<0,-250>[`;\rho]
\square(500,0)|amrb|/->`{@{->}@/^2em/}`->`->/[T`T`Y`Z;1_T`y'`z`S]
\morphism(625,250)|r|<-250,0>[`;\eta]
\morphism(825,375)|r|<0,-250>[`;\sigma]
\efig$$
we have $(\rho\Box\eta)\Box\sigma=\rho\Box(\eta\Box\sigma)$
suggesting that some of the $1_T$-elements of $(SR)(z,x)$ are given by
an obvious coend over $y$ in the category $\bM(T,Y)$.

However, our $\prof$-like notation has its limitations.
For fixed $T$ we can associate to $X$ the category
$\Tilde X=\bM(T,X)$ and to $R\f X\ra Y$ the profunctor
$\Tilde R\f\Tilde X\ra\Tilde Y$ where 
$\Tilde R(y,x)=\bB(T,T)(1_T,y^*Rx)$ but we see no reason why 
a general $1_T$-element of $(SR)(z,x)$ in a general cartesian
bicategory should arise from pasting a $1_T$-element of $S(z,y)$
to a $1_T$-element of $R(y,x)$ for some $y\f T\ra Y$. In short,
while there is a 2-cell $\Tilde S\Tilde R\ra \Tilde{SR}$ in $\prof$
there seems to be no reason why it should have surjective components.
That said, $\Tilde S\Tilde R\ra \Tilde{SR}$ {\em is} an isomorphism
in case $\bB=\spn\E$, for any category $\E$ with finite limits,
and for any cartesian $\bB$ we have isomorphisms 
$1_{\Tilde X}\ra\Tilde{1_X}$ in $\prof$, for any $X$ in $\bB$.
So there is always a normal lax functor
$$\Tilde{(-)}\f\bB\ra\prof$$
which in {\em some} cases is a pseudofunctor. Fortunately, we have
no need for invertibility of the $\Tilde S\Tilde R\ra \Tilde{SR}$.

\subsection{}\label{pl4}
Quite generally, an arrow of $\bG$ as given by the square 
(\ref{Gsquare}) will be called a {\em commutative} square if 
$\alpha$ is invertible.
The arrow (\ref{Gsquare}) of $\bG$ will be said to satisfy the 
{\em Beck-Chevalley condition} if the mate
of $\alpha$ under the adjunctions $f\laj f^*$ and $u\laj u^*$, as given
in the square below (no longer an arrow of $\bG$), is invertible.
$$\bfig
\square(1000,0)/<-`->`->`<-/[X`Y`A`B;f^*`R`S`u^*]
\morphism(1125,250)|a|<250,0>[`;\alpha^*]
\efig$$
Thus Proposition 4.8 of \cite{ckww} says that projection squares of the
form $\tilde p_{R,1_Y}$ and $\tilde r_{1_X,S}$ satisfy the 
Beck-Chevalley condition. (Also, Proposition 4.7 of \cite{ckww} says 
that the same projection squares are commutative. In general,
neither commutative nor Beck-Chevalley implies the other.)
If $R$ and $S$ are also maps and $\alpha$ is invertible then
$\alpha\inv$ gives rise to another arrow of $\bG$ which may or may
not satisfy the Beck-Chevalley condition. The point here is that a
commutative square of maps gives rise to two, generally distinct,
Beck-Chevalley conditions.  It is well known that, for bicategories
of the form $\spn\E$ and $\rel\E$ all pullback
squares of maps satisfy both Beck-Chevalley conditions. A [bi]category
with finite products has automatically a number of pullbacks which
we might call {\em product-absolute} pullbacks because they are
preserved by all [pseudo]functors which preserve products. 

\section{Frobenius Objects in Cartesian Bicategories}\label{Frobenius}
For any object $A$ in $\bB$, we have the following two $\bG$ arrows:

$$\bfig
\square(0,0)/->`->``->/<750,500>%
[A`A\ox A`A\ox A`A\ox(A\ox A);d`d``1\ox d]
\morphism(750,500)|r|<0,-250>[A\ox A`(A\ox A)\ox A;d\ox 1]
\morphism(750,250)|r|<0,-250>[(A\ox A)\ox A`A\ox(A\ox A);a]
\morphism(175,250)|a|<150,0>[`;\_]
\square(1500,0)/->`->`->`/<1500,500>%
[A`A\ox A`A\ox A`A\ox(A\ox A);d`d`1\ox d`]
\morphism(1500,0)|b|<750,0>[A\ox A`(A\ox A)\ox A;d\ox 1]
\morphism(2250,0)|b|<750,0>[(A\ox A)\ox A`A\ox(A\ox A);a]
\morphism(2250,250)|a|<150,0>[`;\_]
\efig$$
obtained from the same equality of arrows in $\map\bB$. (With a
suitable choice of conventions we have equality rather than a mere
isomorphism.) For each square, observe that the data regarded as a 
square in $\bM$ provide an example of a product-absolute 
pullback.

\dfn\label{grpl} An object $A$ is said to be {\em Frobenius\/} if 
both of the $\bG$ arrows above satisfy the Beck-Chevalley condition.
This is to demand invertibility both of 
$\delta_0\f d.d^*\ra 1_A\ox d^*.a.d\ox1_A$, the mate of the first
equality above, and of $\delta_1\f d.d^*\ra d^*\ox1_A.a^*.1_A\ox d$,
the mate of the second equality above.
\eth

\lem\label{sym} The Beck-Chevalley condition for either
square implies the condition for the other.
\eth
\prf Explicitly, in notation suppressing $\ox$, $\delta_0$ and
$\delta_1$ are given by
$$\bfig
\place(-250,250)[\delta_0\; =]
\Vtriangle(0,0)/->`->`<-/[AA`AA`A;1`d^*`d]
\morphism(1000,500)<500,0>[AA`(AA)A;dA]
\morphism(1500,500)<500,0>[(AA)A`A(AA);a]
\morphism(500,0)|b|<1000,0>[A`AA;d]
\Atriangle(1500,0)/<-`->`->/[A(AA)`AA`AA;Ad`Ad^*`1]
\morphism(500,150)|m|<0,225>[`;\epsilon]
\morphism(1250,150)|l|<0,225>[`;|]
\morphism(2000,125)|m|<0,225>[`;A\eta]
\efig$$
and
$$\bfig
\place(-250,250)[\delta_1\; =]
\Vtriangle(0,0)/->`->`<-/[AA`AA`A;1`d^*`d]
\morphism(1000,500)<500,0>[AA`A(AA);Ad]
\morphism(1500,500)<500,0>[A(AA)`(AA)A;a^*]
\morphism(500,0)|b|<1000,0>[A`AA;d]
\Atriangle(1500,0)/<-`->`->/[(AA)A`AA`AA;dA`d^*A`1]
\morphism(500,150)|m|<0,225>[`;\epsilon]
\morphism(1250,150)|l|<0,225>[`;|]
\morphism(2000,125)|m|<0,225>[`;\eta A]
\efig$$

Assume that $\delta_0$ is invertible and paste at its top and right
edges the following pasting composite at its bottom edge.
$$\bfig
\morphism(0,0)|a|<1000,0>[AA`AA;1]
\morphism(1000,0)<500,0>[AA`(AA)A;dA]
\morphism(1500,0)<1500,0>[(AA)A`A(AA);a]
\morphism(0,0)<0,500>[AA`AA;s]
\morphism(1000,0)<0,500>[AA`AA;s]
\morphism(1500,0)<0,500>[(AA)A`A(AA);s]
\morphism(3000,0)<0,500>[A(AA)`(AA)A;s]
\morphism(0,500)|a|<1000,0>[AA`AA;1]
\morphism(1000,500)|a|<500,0>[AA`A(AA);Ad]
\morphism(1500,500)|l|<500,250>[A(AA)`A(AA);As]
\morphism(2000,750)|a|<500,0>[A(AA)`(AA)A;a^*]
\morphism(2500,750)|r|<500,-250>[(AA)A`(AA)A;sA]
\place(500,250)[\iso]
\place(1250,250)[\iso]
\place(2250,250)[\iso]
\morphism(3000,0)<500,0>[A(AA)`AA;Ad^*]
\morphism(3000,500)<500,0>[(AA)A`AA;d^*A]
\morphism(3500,0)<0,500>[AA`AA;s]
\place(3250,250)[\iso]
\efig$$
The squares are pseudonaturality squares for symmetry as 
in 4.5 of \cite{ckww} and the hexagon bounds an invertible
modification constructed from those relating the associativity 
equivalence $a$ and the symmetry equivalence $s$. Next, observe 
that we have $sd\iso d$ and, since $s$ is an equivalence with
$s^*_{A,B}\iso s_{B,A}$, $d^*s\iso d^*$. By functoriality of $\ox$ we
have also $(As)(Ad)\iso Ad$ and $(d^*A)(sA)\iso d^*A$. Noting
the compatibility of the pseudonatural transformation $s$ with the
2-cell $\eta A$, the large pasting composite is seen to be $\delta_1$.
The derivation of invertibility of $\delta_0$ from that of $\delta_1$ 
is effected in a similar way.
\frp

\axm\label{axF}{\em Frobenius}\quad A cartesian bicategory $\bB$
is said to satisfy the {\em Frobenius} axiom if, for each $A$ in $\bB$, 
$A$ is Frobenius.
\eth

\prp\label{gpclofin} In a cartesian bicategory $\bB$,
the Frobenius objects are closed under finite products.
\eth
\prf
Consider a Frobenius object $A$ so that we have invertible
$\delta_0=\delta_0(A)$ in
$$\bfig
\square(0,0)|alrb|/<-`->``<-/<750,500>%
[A`A\ox A`A\ox A`A\ox(A\ox A);d^*`d``1\ox d^*]
\morphism(750,500)|r|<0,-250>[A\ox A`(A\ox A)\ox A;d\ox 1]
\morphism(750,250)|r|<0,-250>[(A\ox A)\ox A`A\ox(A\ox A);a]
\morphism(175,250)|a|<150,0>[`;\delta_0]
\efig$$
For $B$ also Frobenius, form the tensor product of the
diagrams for $\delta_0(A)$ and $\delta_0(B)$, noting that
$\delta_0(A)\ox\delta_0(B)$ is also invertible. The diagram for
$\delta_0(A\ox B)$ is easily formed from that of
$\delta_0(A)\ox\delta_0(B)$
by pasting to its exterior the requisite permutations of
the $A$ and $B$ and using such isomorphisms as
$m(d_A\ox d_B)\iso d_{A\ox B}$, where
$m\f(A\ox A)\ox(B\ox B)\ra (A\ox B)\ox(A\ox B)$ is the middle-four
interchange equivalence.  Thus $A\ox B$ is Frobenius when $A$ and $B$
are so. Invertibility of $\delta_0(I)$ follows easily since $d_I$ is
an equivalence, showing that $I$ is Frobenius.
\frp

Write $\grp\bB$ for the full subbicategory of $\bB$ determined
by the Frobenius objects. It follows immediately from
Proposition \ref{gpclofin} that

\prp\label{eog}
For a cartesian bicategory $\bB$, the full subbicategory
$\grp\bB$ is a cartesian bicategory which satisfies the Frobenius
axiom.
\frp
\eth

In any (pre)cartesian bicategory we have, for each object $X$, the 
following arrows: 
$$N_X\quad=\quad I\to^{t_X^*}X\to^{d_X}X\ox X\quad\quad\mbox{and}%
\quad\quad E_X\quad=\quad X\ox X\to^{d_X^*}X\to^{t_X}I$$ 
Since the cartesian bicategory $\bB$ is a (symmetric) monoidal
bicategory it can be seen as a one-object tricategory,
so that pseudo adjunctions $N,E\f X\laj A$, where $X$ and $A$ are
{\em objects} of $\bB$ (and $N$ and $E$ are arrows of $\bB$), are 
well defined. (We note that, especially since $\bB$ is symmetric, it
is customary to speak of such $X$ and $A$ as {\em duals}.)

\prp\label{XlajX}
For a Frobenius object $X$ in a cartesian bicategory,
$N_X$ and $E_X$ provide the unit and counit for a pseudo-adjunction
$X\laj X$.  
\eth
\prf (Sketch)
We are to exhibit isomorphisms
$$(E_X\ox X)a^*(X\ox N_X)\iso s_{X,I}\quad\mbox{and}\quad%
(X\ox E_X)a(N_X\ox X)\iso s_{I,X}$$
subject to two coherence equations.
Consider:
$$\bfig
\morphism(0,0)|a|<500,0>[X\ox I`X\ox X;X\ox t_X^*]
\morphism(500,0)|a|<750,0>[X\ox X`X\ox(X\ox X);X\ox d_X]
\morphism(1250,0)|a|<750,0>[X\ox(X\ox X)`(X\ox X)\ox X;a^*]
\morphism(2000,0)|r|<0,-500>[(X\ox X)\ox X`X\ox X;d_X^*\ox X]
\morphism(500,0)|r|<0,-500>[X\ox X`X;d_X^*]
\morphism(500,-500)|b|<1500,0>[X`X\ox X;d_X]
\morphism(0,0)|l|<500,-500>[X\ox I`X;r]
\morphism(2000,-500)|r|<0,-500>[X\ox X`I\ox X;t_X\ox X]
\morphism(500,-500)|b|<1500,-500>[X`I\ox X;l]
\morphism(1100,-300)|a|<200,0>[`;\stackrel{\delta_1}{\simeq}]
\place(350,-200)[\iso]
\place(1750,-750)[\iso]
\efig$$
$$\bfig
\morphism(0,0)|a|<500,0>[I\ox X`X\ox X;t_X^*\ox X]
\morphism(500,0)|a|<750,0>[X\ox X`(X\ox X)\ox X;d_X\ox X]
\morphism(1250,0)|a|<750,0>[(X\ox X)\ox X`X\ox(X\ox X);a]
\morphism(2000,0)|r|<0,-500>[X\ox(X\ox X)`X\ox X;X\ox d_X^*]
\morphism(500,0)|r|<0,-500>[X\ox X`X;d_X^*]
\morphism(500,-500)|b|<1500,0>[X`X\ox X;d_X]
\morphism(0,0)|l|<500,-500>[I\ox X`X;l]
\morphism(2000,-500)|r|<0,-500>[X\ox X`X\ox I;X\ox t_X]
\morphism(500,-500)|b|<1500,-500>[X`X\ox I;r]
\morphism(1300,-300)|a|<-200,0>[`;\stackrel{{\delta_0}\inv}{\simeq}]
\place(350,-200)[\iso]
\place(1750,-750)[\iso]
\efig$$
For the coherence requirements let us abbreviate $\ox$ by
juxtaposition, as we have before, but now work as if the bicategory 
constraints of $\bB$ and those of the monoidal structure 
$(\bB,\ox,I)$ are strict. (In general, this is not acceptable
because a monoidal bicategory is not tri-equivalent to a one-object
3-category. However, our monoidal structure, being given by universal
properties, is less problematical.) Temporarily, write
$N\f I\ra X^\circ X$ and $E\f XX^\circ\ra I$, just to mark
the role of the $X$'s. Write $\alpha\f 1_X\ra (EX)(XN)$ and 
$\beta\f (X^\circ E)(NX^\circ)\ra1_{X^\circ}$
for the isomorphisms built from those above, with the
simplifying assumptions. The coherence requirements of $\alpha$ 
and $\beta$ are the pasting equations
$$\bfig
\place(-1100,1000)[1_E =]
\square(0,0)<600,600>[XX^\circ XX^\circ`XX^\circ`XX^\circ`I;%
XX^\circ E`EXX^\circ`E`E]
\morphism(-600,1200)|m|<600,-600>[XX^\circ`XX^\circ XX^\circ;XNX^\circ]
\morphism(-600,1200)|l|/{@{->}@/_3em/}/<600,-1200>%
[XX^\circ`XX^\circ;XX^\circ]
\morphism(-600,1200)|r|/{@{->}@/^3em/}/<1200,-600>%
[XX^\circ`XX^\circ;XX^\circ]
\morphism(-500,600)|a|<200,0>[`;\alpha X^\circ]
\morphism(0,900)|a|<200,0>[`;X\beta]
\place(300,300)[\iso]
\square(1200,600)<600,600>[I`X^\circ X`X^\circ X`X^\circ XX^\circ X;%
N`N`NX^\circ X`X^\circ XN]
\place(1500,900)[\iso]
\morphism(1800,1200)|r|/{@{->}@/^3em/}/<600,-1200>%
[X^\circ X`X^\circ X;X^\circ X]
\morphism(1200,600)|l|/{@{->}@/_3em/}/<1200,-600>%
[X^\circ X`X^\circ X;X^\circ X]
\morphism(1800,600)|m|<600,-600>%
[X^\circ XX^\circ X`X^\circ X;X^\circ EX] 
\morphism(1500,300)|a|<200,0>[`;X^\circ\alpha]
\morphism(2100,600)|a|<200,0>[`;\beta X]
\place(700,1000)[1_N =]
\efig$$
where the unlabelled isomorphisms in the squares are given by
pseudofunctoriality of $\ox$. 
We will verify the first of these equations, verification of the 
second being similar, now using $X^\circ = X$ but continuing to
supress the constraints both for $\bB$ and for the monoidal structure.
Thus we must show that the composite on the left below
$$\bfig
\square(0,0)|ammb|[XXXX`XXX`XXX`XX;XXd^*`d^*XX`d^*X`Xd^*]
\place(250,250)[\iso]
\morphism(-350,850)|m|<350,-350>[XXX`XXXX;XdX]
\morphism(-350,850)|l|<0,-500>[XXX`XX;d^*X]
\morphism(-350,850)|a|<500,0>[XXX`XX;Xd^*]
\morphism(150,850)|r|<350,-350>[XX`XXX;Xd]
\morphism(-350,350)|b|<350,-350>[XX`XXX;dX]
\morphism(-225,350)|a|<150,0>[`;\delta_1X]
\morphism(0,675)|a|<150,0>[`;X\delta_0\inv]
\morphism(-700,1200)|m|<350,-350>[XX`XXX;Xt^*X]
\morphism(500,0)|m|<350,-350>[XX`I;tt]
\place(1075,500)[=]
\morphism(1300,1200)|m|<350,-350>[XX`XXX;Xt^*X]
\square(1650,350)|ammb|[XXX`XX`XX`X;Xd^*`d^*X`d^*`d^*]
\place(1900,600)[\iso]
\morphism(2150,350)|m|<350,-350>[X`XX;d]
\morphism(1650,350)|b|<350,-350>[XX`XXX;dX]
\morphism(2150,850)|r|<350,-350>[XX`XXX;Xd]
\morphism(2000,0)|b|<500,0>[XXX`XX;Xd^*]
\morphism(2500,500)|r|<0,-500>[XXX`XX;d^*X]
\morphism(2500,0)|m|<350,-350>[XX`I;tt]
\morphism(2250,350)|a|<150,0>[`;\delta_1]
\morphism(2000,175)|a|<150,0>[`;\delta_0\inv]
\efig$$
is $1_E$. Again using pseudofunctoriality of $\ox$, we have 
the equality shown and finally the diagram on the right can 
be shown to be $1_E$ from the definitions of $\delta_0$ and $\delta_1$.
\frp

\subsection{}\label{universal}
If $R\f X\ra A$ is an arrow in $\bB$ then given
pseudo adjunctions $X\laj X^\circ$ and $A\laj A^\circ$  
we should expect that adaption of the calculus of 
{\em mates} found in \cite{kas} will enable us to define 
$R^\circ\f X^\circ\ra A^\circ$ by the usual formula. In fact, if
every object of $\bB$ has a dual one should expect 
$(-)^{\circ}$ to provide a pseudofunctor 
$(-)^{\circ}\f\bB\oprev\ra\bB$ between tricategories, where $(-)\rev$
denotes dualization with respect to objects of $\bB$ {\em composed}
via $\ox$, while as usual  $(-)\op$ denotes dualization with respect 
to the 1-cells of $\bB$. In
particular, one should expect 
$(X\ox Y)^\circ\simeq Y^\circ\ox X^\circ$. 
The point of this paragraph is that the $(-)^{\circ}$ of the following
proposition arises from the properties already under consideration
and is not a new structure as in the similarly denoted operation
of \cite{fas}.

\prp\label{ring}
For a cartesian bicategory $\bB$ in which every object is
Frobenius, there is an involutory pseudofunctor 
$$(-)^{\circ}\f\bB\op\ra\bB$$
which is the identity on objects.
\eth
\prf
With $X^\circ = X$
we define 
$$(-)^\circ_{A,X}\f\bB\op(A,X)=\bB(X,A)\ra\bB(A,X)$$
by the evidently
functorial formula 
$$R^\circ = (X\ox E_A)(X\ox R\ox A)(N_X\ox A)$$
In terms of the one object tricategory $(\bB,\ox,I)$ with
single object $*$, we can express $R^\circ$ by the pasting
$$\bfig
\scalefactor{0.75}
\Vtriangle(1500,0)|tlr|/->`->`<-/[*`*`*;I`A`A]
\morphism(2000,150)|m|<0,200>[`;E_A]
\morphism(2000,0)|b|<1000,0>[*`*;I]
\morphism(3000,0)|r|<500,500>[*`*;X]
\morphism(2500,500)|a|<1000,0>[*`*;I]
\morphism(2750,150)|m|<0,200>[`;R]
\Atriangle(3000,0)|rlb|/<-`->`->/[*`*`*;X`X`I]
\morphism(3500,150)|m|<0,200>[`;N_X]
\efig$$
For $R\f X\ra A$, along with $S\f A\ra Y$, to give 
$\widetilde{(-)^{\circ}}\f R^\circ S^\circ\ra (SR)^\circ$ we consider
$$\bfig
\scalefactor{0.75}
\morphism(-500,500)|l|<500,-500>[*`*;Y]
\morphism(0,0)|r|<500,500>[*`*;Y]
\morphism(0,0)|b|<1000,0>[*`*;I]
\morphism(0,150)|m|<0,200>[`;E_Y]
\morphism(-500,500)|a|<1000,0>[*`*;I]
\morphism(1000,0)|r|<500,500>[*`*;]
\morphism(500,500)|a|<1000,0>[*`*;I]
\morphism(750,150)|m|<0,200>[`;S]
\Atriangle(1000,0)|rlb|/<-`->`->/[*`*`*;A`A`I]
\morphism(1500,150)|m|<0,200>[`;N_A]
\Vtriangle(1500,0)/->`->`<-/[*`*`*;I``A]
\morphism(2000,150)|m|<0,200>[`;E_A]
\morphism(2000,0)|b|<1000,0>[*`*;I]
\morphism(3000,0)|r|<500,500>[*`*;X]
\morphism(2500,500)|a|<1000,0>[*`*;I]
\morphism(2750,150)|m|<0,200>[`;R]
\Atriangle(3000,0)|rlb|/<-`->`->/[*`*`*;X`X`I]
\morphism(3500,150)|m|<0,200>[`;N_X]
\efig$$
in which the pasting composite displays $R^\circ S^\circ$.
The required $\widetilde{(-)^{\circ}}$ is obtained as the
collapsing of the centre triangles using
$\alpha\inv\f(E_A\ox A)(A\ox N_A)\iso s_{A,I}$ of the
pseudo adjunction $N_A,E_A\f A\laj A$. Evidently,
$\widetilde{(-)^{\circ}}$ is invertible. We give the identity
constraint for $(-)^{\circ}$ as
$\beta\inv\f1_X\ra(X\ox E_X)(N_X\ox X)$ which is again invertible.
Finally, having observed that the mate
description of $R^\circ=(X\ox E_A)(X\ox R\ox A)(N_X\ox A)$ was
given by expanding $R\f X\ra A$ as $R\f X\ox I\ra I\ox A$ we see
by writing $R\f I\ox X\ra A\ox I$ that we have equally
$$R^\circ \iso\ (E_A\ox X)(A\ox R\ox X)(A\ox N_X)$$
Thus we may as well give $$((-)^{\circ})\op\f\bB\ra\bB\op$$
by the formula 
$$(A\to^S X)\rass (E_X\ox A)(X\ox S\ox A)(X\ox N_A)$$
so that $R^{\circ\circ}$ is the pasting
$$\bfig
\scalefactor{0.75}
\Atriangle(1000,0)/<-`->`->/[*`*`*;A`A`I]
\morphism(1500,150)|m|<0,200>[`;N_A]
\Vtriangle(1500,0)/->`->`<-/[*`*`*;I``A]
\morphism(2000,150)|m|<0,200>[`;E_A]
\morphism(2000,0)|b|<1000,0>[*`*;I]
\morphism(3000,0)|l|<500,500>[*`*;X]
\morphism(2500,500)|a|<1000,0>[*`*;I]
\morphism(2750,150)|m|<0,200>[`;R]
\Atriangle(3000,0)/<-`->`->/[*`*`*;X``I]
\morphism(3500,150)|m|<0,200>[`;N_X]
\Vtriangle(3500,0)/->`->`<-/[*`*`*;I`X`X]
\morphism(4000,150)|m|<0,200>[`;E_X]
\efig$$
and we have a canonical isomorphism $R\iso R^{\circ\circ}$,
again using the $\alpha$ and $\beta$ constraints of the pseudo
adjunctions $N_X,E_X\f X\laj X$ of Proposition \ref{XlajX}.
\frp

\prp\label{N_RandE_R}
For an arrow $R\f X\ra A$ in a cartesian bicategory, with
$X$ and $A$ Frobenius, if the $\tilde d_R$ and $\tilde t_R$ of the units
$$\bfig
\square(0,0)[X`X\ox X`A`A\ox A;d_X`R`R\ox R`d_A]
\morphism(125,250)|m|<250,0>[`;\tilde d_R]
\square(1250,0)[X`I`A`I;t_X`R`\top`t_A]
\morphism(1375,250)|m|<250,0>[`;\tilde t_R]
\efig$$
are invertible then we can construct squares $N_R$ and $E_R$
$$\bfig
\square(0,0)|alrm|[I`I`X`A;1_I`t^*_X`t^*_A`R]
\morphism(125,250)|m|<250,0>[`;\tilde t^*_R]
\square(0,-500)|mlrb|[X`A`X\ox X`A\ox A;R`d_X`d_A`R\ox R]
\morphism(125,-250)|m|<250,0>[`;\tilde d\inv_R]
\morphism(0,500)|l|/{@{->}@/_2em/}/<0,-1000>[I`X\ox X;N_X]
\morphism(500,500)|m|/{@{->}@/^2em/}/<0,-1000>[I`A\ox A;N_A]
\place(-500,250)[N_R =]
\square(1500,0)|alrm|[X\ox X`A\ox A`X`A;R\ox R`d^*_X`d^*_A`R]
\morphism(1625,250)|m|<250,0>[`;\tilde d^*_R]
\square(1500,-500)|mlrb|[X`A`I`I;R`t_X`t_A`1_I]
\morphism(1625,-250)|m|<250,0>[`;\tilde t_R\inv]
\morphism(1500,500)|m|/{@{->}@/_2em/}/<0,-1000>[X\ox X`I;E_X]
\morphism(2000,500)|r|/{@{->}@/^2em/}/<0,-1000>[A\ox A`I;E_A]
\place(1000,250)[E_R =]
\square(2500,-250)|alrm|[X`A`X`A;R`X`A`R]
\morphism(2625,0)|m|<250,0>[`;R]
\efig$$
where $\tilde t^*_R$ is the mate of $\tilde t_R$ and 
$\tilde d^*_R$ is the mate of $\tilde d_R$,
which when tensored with the identity square $R$, above, satisfy
the following equations (in which $\ox$ is suppressed):
\begin{equation}\label{NRR/RER}
\bfig
\square(0,0)|alra|[X`A`XXX`AAA;R`N_XX`N_AA`RRR]
\morphism(125,250)|a|<250,0>[`;N_RR]
\square(0,-500)|alrb|[XXX`AAA`X`A;RRR`XE_X`AE_A`R]
\morphism(125,-250)|a|<250,0>[`;RE_R]
\morphism(0,500)|l|/{@{->}@/_3em/}/<0,-1000>[X`X;X]
\place(-200,0)[\iso]
\morphism(500,500)|r|/{@{->}@/^3em/}/<0,-1000>[A`A;A]
\place(700,0)[\iso]
\square(1500,0)|alra|[X`A`XXX`AAA;R`XN_X`AN_A`RRR]
\morphism(1625,250)|a|<250,0>[`;RN_R]
\square(1500,-500)|alrb|[XXX`AAA`X`A;RRR`E_XX`E_AA`R]
\morphism(1625,-250)|a|<250,0>[`;E_RR]
\morphism(1500,500)|l|/{@{->}@/_3em/}/<0,-1000>[X`X;X]
\place(1300,0)[\iso]
\morphism(2000,500)|r|/{@{->}@/^3em/}/<0,-1000>[A`A;A]
\place(2200,0)[\iso]
\place(-500,250)[R =]
\place(2500,250)[= R]
\efig
\end{equation}
\eth
\prf
The vertical edges of the diagrams have been clarified
in Proposition \ref{XlajX}. For the rest it suffices for 
each equation to expand $N_R$ and $E_R$, verify the following
equalities
$$\bfig
\Atriangle(0,0)|lrl|/->`->`/<250,250>[XX`X`XXX;d^*`dX`]
\Vtriangle(0,-250)|llr|/`->`->/<250,250>[X`XXX`XX;`d`Xd^*]
\morphism(150,0)|m|<150,0>[`;\delta_0]
\morphism(250,250)|a|<750,0>[XX`AA;RR]
\morphism(650,125)|m|<250,0>[`;\tilde d\inv_RR]
\morphism(500,0)|m|<750,0>[XXX`AAA;RRR]
\morphism(650,-125)|m|<250,0>[`;R\tilde d^*_R]
\morphism(250,-250)|b|<750,0>[XX`AA;RR]
\morphism(1000,250)|r|<250,-250>[AA`AAA;dA]
\morphism(1250,0)|r|<-250,-250>[AAA`AA;Ad^*]
\place(1550,0)[=]
\morphism(2000,250)|a|<750,0>[XX`AA;RR]
\morphism(2125,125)|m|<250,0>[`;\tilde d^*_R]
\morphism(1750,0)|m|<750,0>[X`A;R]
\morphism(2125,-125)|m|<250,0>[`;\tilde d\inv_R]
\morphism(2000,-250)|a|<750,0>[XX`AA;RR]
\morphism(2000,250)|l|<-250,-250>[XX`X;d^*]
\morphism(2750,250)|a|<250,-250>[AA`AAA;dA]
\morphism(3000,0)|r|<-250,-250>[AAA`AA;Ad^*]
\morphism(2750,250)|a|<-250,-250>[AA`A;d^*]
\morphism(2675,0)|m|<150,0>[`;\delta_0]
\morphism(2500,0)|l|<250,-250>[A`AA;d]
\morphism(1750,0)|l|<250,-250>[X`XX;d]
\efig$$
$$\bfig
\Atriangle(0,0)|lrl|/->`->`/<250,250>[XX`X`XXX;d^*`Xd`]
\Vtriangle(0,-250)|llr|/`->`->/<250,250>[X`XXX`XX;`d`d^*X]
\morphism(150,0)|m|<150,0>[`;\delta_1]
\morphism(250,250)|a|<750,0>[XX`AA;RR]
\morphism(650,125)|m|<250,0>[`;R\tilde d\inv_R]
\morphism(500,0)|m|<750,0>[XXX`AAA;RRR]
\morphism(650,-125)|m|<250,0>[`;\tilde d^*_RR]
\morphism(250,-250)|b|<750,0>[XX`AA;RR]
\morphism(1000,250)|r|<250,-250>[AA`AAA;Ad]
\morphism(1250,0)|r|<-250,-250>[AAA`AA;d^*A]
\place(1550,0)[=]
\morphism(2000,250)|a|<750,0>[XX`AA;RR]
\morphism(2125,125)|m|<250,0>[`;\tilde d^*_R]
\morphism(1750,0)|m|<750,0>[X`A;R]
\morphism(2125,-125)|m|<250,0>[`;\tilde d\inv_R]
\morphism(2000,-250)|a|<750,0>[XX`AA;RR]
\morphism(2000,250)|l|<-250,-250>[XX`X;d^*]
\morphism(2750,250)|a|<250,-250>[AA`AAA;Ad]
\morphism(3000,0)|r|<-250,-250>[AAA`AA;d^*A]
\morphism(2750,250)|a|<-250,-250>[AA`A;d^*]
\morphism(2675,0)|m|<150,0>[`;\delta_1]
\morphism(2500,0)|l|<250,-250>[A`AA;d]
\morphism(1750,0)|l|<250,-250>[X`XX;d]
\efig$$
and use such further equalities as
$$\bfig
\square(0,0)|alra|[X`A`XX`AA;R`d`d`RR]
\morphism(125,250)|a|<250,0>[`;\tilde d_R\inv]
\square(0,-500)|alrb|[XX`AA`X`A;RR`Xt`At`R]
\morphism(125,-250)|a|<250,0>[`;R\tilde t_R\inv]
\morphism(0,500)|l|/{@{->}@/_3em/}/<0,-1000>[X`X;X]
\place(-200,0)[\iso]
\morphism(500,500)|r|/{@{->}@/^3em/}/<0,-1000>[A`A;A]
\place(700,0)[\iso]
\square(1500,0)|alra|[X`A`XX`AA;R`t^*X`t^*A`RR]
\morphism(1625,250)|a|<250,0>[`;\tilde t^*_RR]
\square(1500,-500)|alrb|[XX`AA`X`A;RR`d^*`d^*`R]
\morphism(1625,-250)|a|<250,0>[`;\tilde d_R^*]
\morphism(1500,500)|l|/{@{->}@/_3em/}/<0,-1000>[X`X;X]
\place(1300,0)[\iso]
\morphism(2000,500)|r|/{@{->}@/^3em/}/<0,-1000>[A`A;A]
\place(2200,0)[\iso]
\place(-500,250)[R =]
\place(1000,250)[R =]
\efig$$
\frp
\subsection{}\label{comonoids}
Every object $X$ of a bicategory with finite products is, 
essentially uniquely, a pseudo comonoid via $d_X$ and $t_X$. 
It follows that every object $X$ in a cartesian bicategory $\bB$
is a (pseudo) comonoid (via $d_X$ and $t_X$) since $\bM$ has
finite products and the inclusion functor $i\f\bM\ra\bB$ 
is strongly monoidal. (It is the identity on objects and we 
observe from Proposition 3.24 of \cite{ckww} that 
$f\x g{\s\to^\simeq}f\ox g$ in $\bB$.)
Similarly, for $R\f X\ra A$ in $\bB$,
$R$ has an essentially unique comonoid structure in $\bG$, via
$(d_X,\tilde d_R,d_A)$ and $(t_X,\tilde t_R,t_A)$, 
since $\bG$ has finite products. In fact, given $d_X$ and $d_A$,
$\tilde d_R$ is uniquely determined and given $t_X$ and $t_A$,
$\tilde t_R$ is uniquely determined. This fact can be reinterpretted 
to say that $R\f X\ra A$ has an essentially unique
lax comonoid homomorphism structure via $d_R=(d_X,\tilde d_R,d_A)$
and $t_R=(t_X,\tilde t_R,t_A)$ which is then a 
{\em comonoid homomorphism} if and only if the 2-cells $\tilde d_R$ and 
$\tilde t_R$ are invertible. Thus being a comonoid homomorphism is a 
{\em property} of an arrow in a cartesian bicategory.

\thm\label{mapiff}
For an arrow $R\f X\ra A$ in a cartesian bicategory, with
$X$ and $A$ Frobenius, the following are equivalent:
\begin{enumerate}
\item[(1)] $R$ is a map;
\item[(2)] $R$ is a comonoid homomorphism;
\item[(3)] $R\laj R^\circ$.
\end{enumerate}
\eth
\prf
(1) implies (2) follows from the fact that $d$ and $t$ are 
pseudonatural on maps and (3) implies (1) is trivial.
So, assuming (2), that $R$ is a comonoid homomorphism,
construct $N_R$ and $E_R$ as in Proposition \ref{N_RandE_R} and
define (suppressing $\ox$ as usual)
$$\bfig
\place(-750,0)[\eta_R=]
\morphism(0,500)|l|<0,-500>[X`XXX;N_XX]
\morphism(0,0)|m|<350,350>[XXX`XXA;XXR]
\morphism(0,0)|m|<350,-350>[XXX`XAA;XRR]
\morphism(0,0)|l|<0,-500>[XXX`X;XE_X]
\morphism(350,850)|r|<0,-500>[A`XXA;N_XA]
\morphism(350,350)|r|<0,-700>[XXA`XAA;XRA]
\morphism(350,-350)|r|<0,-500>[XAA`X;XE_A]
\morphism(0,500)|a|<350,350>[X`A;R]
\morphism(0,-500)|b|<350,-350>[X`X;X]
\morphism(150,550)|b|<150,0>[`;N_XR]
\place(250,0)[\iso]
\morphism(150,-550)|a|<150,0>[`;XE_R]
\morphism(0,500)|l|/{@{->}@/_3em/}/<0,-1000>[X`X;X]
\place(-200,0)[\iso]
\place(1000,0)[\epsilon_R=]
\morphism(1500,850)|l|<0,-500>[A`XXA;N_XA]
\morphism(1500,350)|l|<0,-700>[XXA`XAA;XRA]
\morphism(1500,-350)|l|<0,-500>[XAA`X;XE_A]
\morphism(1500,850)|a|<350,-350>[A`A;A]
\morphism(1500,350)|m|<350,-350>[XXA`AAA;RRA]
\morphism(1500,-350)|m|<350,350>[XAA`AAA;RAA]
\morphism(1500,-850)|r|<350,350>[X`A;R]
\morphism(1850,500)|r|<0,-500>[A`AAA;N_AA]
\morphism(1850,0)|r|<0,-500>[AAA`A;AE_A]
\morphism(1850,500)|r|/{@{->}@/^3em/}/<0,-1000>[A`A;A]
\morphism(1600,550)|b|<150,0>[`;N_RA]
\place(1600,0)[\iso]
\place(2050,0)[\iso]
\morphism(1600,-550)|a|<150,0>[`;RE_A]
\efig$$
where we note that both three-fold vertical composites
are the arrow $R^\circ$, $N_XR=1_{N_X}\ox1_R$ and
$RE_A=1_R\ox1_{E_A}$ are isomorphisms while
$XE_R=1_{1_X}\ox E_R$ and $N_RA=N_R\ox1_{1_A}$.
When $\eta_R$ and $\epsilon_R$ are pasted at $R^\circ$
the result is
$$\bfig
\morphism(0,500)|l|<0,-500>[X`XXX;N_XX]
\morphism(0,0)|m|<350,350>[XXX`XXA;XXR]
\morphism(0,0)|m|<350,-350>[XXX`XAA;XRR]
\morphism(0,0)|m|<0,-500>[XXX`X;XE_X]
\morphism(350,850)|m|<0,-500>[A`XXA;N_XA]
\morphism(350,-350)|m|<0,-500>[XAA`X;XE_A]
\morphism(0,500)|a|<350,350>[X`A;R]
\morphism(0,-500)|b|<350,-350>[X`X;X]
\morphism(150,550)|b|<150,0>[`;N_XR]
\place(350,0)[\iso]
\morphism(150,-550)|a|<150,0>[`;XE_R]
\morphism(0,500)|l|/{@{->}@/_3em/}/<0,-1000>[X`X;X]
\place(-200,0)[\iso]
\morphism(350,850)|a|<350,-350>[A`A;A]
\morphism(350,350)|m|<350,-350>[XXA`AAA;RRA]
\morphism(350,-350)|m|<350,350>[XAA`AAA;RAA]
\morphism(350,-850)|r|<350,350>[X`A;R]
\morphism(700,500)|r|<0,-500>[A`AAA;N_AA]
\morphism(700,0)|r|<0,-500>[AAA`A;AE_A]
\morphism(700,500)|r|/{@{->}@/^3em/}/<0,-1000>[A`A;A]
\morphism(450,550)|b|<150,0>[`;N_RA]
\morphism(450,-550)|a|<150,0>[`;RE_A]
\place(900,0)[\iso]
\place(1350,0)[=]
\square(2000,0)|alra|[X`A`XXX`AAA;R`N_XX`N_AA`RRR]
\morphism(2125,250)|a|<250,0>[`;N_RR]
\square(2000,-500)|alrb|[XXX`AAA`X`A;RRR`XE_X`AE_A`R]
\morphism(2125,-250)|a|<250,0>[`;RE_R]
\morphism(2000,500)|l|/{@{->}@/_3em/}/<0,-1000>[X`X;X]
\place(1800,0)[\iso]
\morphism(2500,500)|r|/{@{->}@/^3em/}/<0,-1000>[A`A;A]
\place(2700,0)[\iso]
\place(3000,0)[=]
\place(3200,0)[R]
\efig$$
the first equality from functoriality of $\ox$, the 
second equality being the first equation of (\ref{NRR/RER})
of Proposition \ref{N_RandE_R}.
To complete the proof that we have an adjunction 
$\eta_R,\epsilon_R\f R\laj R^\circ$ we must show that
when $\eta_R$ is pasted to $\epsilon_R$ at $R$ the result
is $R^\circ$. To aid readability we draw as commutative as
many regions as possible. Consider:
$$\bfig
\morphism(0,0)|b|<500,0>[A`XXA;N_XA]
\morphism(0,0)<250,500>[A`A;A]
\morphism(250,500)<500,0>[A`AAA;N_AA]
\morphism(300,125)|m|<0,250>[`;N_RA]
\Atriangle(500,0)/<-`<-`->/<250,500>[AAA`XXA`XAA;RRA`RAA`XRA]
\morphism(1000,0)|b|<500,0>[XAA`X;XE_A]
\morphism(750,500)<500,0>[AAA`A;AE_A]
\morphism(1500,0)|m|<-250,500>[X`A;R]
\morphism(0,0)|b|/{@{->}@/_3em/}/<1500,0>[A`X;R^\circ]
\morphism(250,500)/{@{->}@/^3em/}/<1000,0>[A`A;A]
\morphism(1500,0)|b|<500,0>[X`XXX;N_XX]
\morphism(1250,500)<500,0>[A`XXA;N_XA]
\Vtriangle(1750,0)/->`<-`<-/<250,500>[XXA`XAA`XXX;XRA`XXR`XRR]
\morphism(2000,0)|b|<500,0>[XXX`X;XE_X]
\morphism(2250,500)<500,0>[XAA`X;XE_A]
\morphism(2500,0)|r|<250,500>[X`X;X]
\morphism(2450,125)|m|<0,250>[`;XE_R]
\morphism(1500,0)|b|/{@{->}@/_3em/}/<1000,0>[X`X;X]
\morphism(1250,500)/{@{->}@/^3em/}/<1500,0>[A`X;R^\circ]
\efig$$
(which is the requisite pasting rotated 90 degrees counterclockwise). 
Rearrange it as below:
\begin{equation}\label{pasting}
\bfig
\square(0,0)|alrb|<500,500>[A`A`XXA`AAA;A`N_XA`N_AA`RRA]
\morphism(125,250)|a|<250,0>[`;N_RA]
\square(500,0)|armb|<750,500>[A`XXA`AAA`XXAAA;N_XA`N_AA`XXN_AA`N_XAAA]
\square(1250,0)|ammb|<750,500>%
[XXA`XAA`XXAAA`XAAAA;XRA`XXN_AA`XAN_AA`XRAAA]
\square(1250,-500)|bmma|/->`<-`<-`->/<750,500>%
[XXAAA`XAAAA`XXXXA`XXXAA;XRAAA`XXRRA`XRRAA`XXXRA]
\square(1250,-1000)|ammb|<750,500>%
[XXXXA`XXXAA`XXA`XAA;XXXRA`XE_XXA`XE_XAA`XRA]
\square(2000,-1000)|amlb|<750,500>%
[XXXAA`XXX`XAA`X;XXXE_A`XE_XAA`XE_X`XE_A]
\square(2750,-1000)<500,500>[XXX`XAA`X`X;XRR`XE_X`XE_A`X]
\morphism(2875,-750)|a|<250,0>[`;XE_R]
\morphism(1250,-500)|m|<750,500>[XXXXA`XAAAA;XRRRA]
\morphism(0,0)|m|<1250,-500>[XXA`XXXXA;N_XXXA]
\morphism(0,0)|m|<1250,-1000>[XXA`XXA;XXA]
\morphism(2000,0)|m|<1250,-500>[XAAAA`XAA;XAAE_A]
\morphism(2000,500)|m|<1250,-1000>[XAA`XAA;XAA]
\efig
\end{equation}
The following prism commutes:
$$\bfig
\square(0,0)<750,500>[A`A`XXA`AAA;A`N_XA`N_AA`RRA]
\morphism(250,250)|a|<250,0>[`;N_RA]
\square(1500,-500)<750,500>[XXA`XXA`XXXXA`XXAAA;XXA`XXN_XA`XXN_AA`XXRRA]
\morphism(1750,-250)|a|<250,0>[`;XXN_RA]
\morphism(0,500)||/@{->}_<>(0.8){N_XA}/<1500,-500>[A`XXA;]
\morphism(750,500)||/@{->}^<>(0.7){N_XA}/<1500,-500>[A`XXA;]
\morphism(0,0)||/@{->}_<>(0.3){N_XXXA}/<1500,-500>[XXA`XXXXA;]
\morphism(750,0)||/@{->}^<>(0.2){N_XAAA}/<1500,-500>[AAA`XXAAA;]
\morphism(3000,-500)|r|<0,-500>[XAA`XAAAA;XAN_AA]
\morphism(2250,0)||/@{->}^<>(0.8){XRA}/<750,-500>[XXA`XAA;]
\morphism(1500,0)||/@{->}_<>(0.8){XRA}/<1500,-500>[XXA`XAA;]
\morphism(2250,-500)||/@{->}^<>(0.3){XRAAA}/<750,-500>[XXAAA`XAAAA;]
\morphism(1500,-500)||/@{->}_<>(0.5){XRRRA}/<1500,-500>[XXXXA`XAAAA;]
\morphism(2700,-750)|a|<250,0>[`;XRN_RA]
\efig$$
Replace the top three squares of (\ref{pasting}) above by the two
front faces of the prism. Employ a similar commuting prism
to replace the bottom three squares of (\ref{pasting}) and
obtain:
$$\bfig
\square(-750,0)|alla|<750,500>[A`XXA`XXA`XXXXA;N_XA`N_XA``N_XXXA]
\morphism(-750,0)|l|<750,-500>[XXA`XXA;XXA]
\square(0,0)|alra|<750,500>[XXA`XAA`XXXXA`XAAAA;XRA`XXN_XA`XAN_AA`XRRRA]
\morphism(250,250)|a|<250,0>[`;XRN_RA]
\square(0,-500)|alrb|<750,500>[XXXXA`XAAAA`XXA`XAA;%
                               XRRRA`XE_XXA`XE_AAA`XRA]
\morphism(250,-250)|a|<250,0>[`;XE_RRA]
\square(750,-500)|alrb|<750,500>[XAAAA`XAA`XAA`X;XAAE_A``XE_A`XE_A]
\morphism(750,500)<750,-500>[XAA`XAA;XAA]
\efig$$
$$\bfig
\place(-500,250)[=]
\square(0,0)|allb|[A`XXA`A`XXA;N_XA`A`XXA`N_XA]
\square(500,0)|alrb|[XXA`XAA`XXA`XAA;XRA`XXA`XAA`XRA]
\square(1000,0)|arrb|[XAA`X`XAA`X;XE_A`XAA`X`XE_A]
\place(2000,250)[= R^\circ]
\efig$$
where the penultimate equality is obtained from the second 
equation of (\ref{NRR/RER}) of Proposition~\ref{N_RandE_R}
by tensoring it on the left by $X$
and on the right by $A$ and applying the result to the two
middle squares of the penultimate pasting.
\frp

\subsection{}\label{*=circ}
From Theorem \ref{mapiff} it follows that for a 
map $f\f X\ra A$, with $X$ and $A$ Frobenius in a cartesian 
bicategory, we have $f^*\iso f^\circ$ and we may as well write
$f^*=f^\circ$ for our specified right adjoints in this event
and use the explicit formula for $f^\circ$ when it is convenient
to do so. 

\thm\label{grplhom}
If $A$ is a Frobenius object in a cartesian bicategory $\bB$,
then, for all $T$ in $\bB$, the hom-category $\bM(T,A)$ is a groupoid.
\eth

We will break the proof of Theorem \ref{grplhom} into a sequence
of lemmas and employ the notation of \ref{pl3}. 

\lem\label{rewritdelta}
With reference to the 2-cell $\delta_1$ in Definition \ref{grpl},
$$dd^*{\s {\to^\simeq}}(p^*\wedge r^*)(p\wedge r)\quad\mbox{and}\quad%
(d^*\ox X)(X\ox d){\s {\to^\simeq}}p^*p\wedge p^*r\wedge r^*r$$
and these canonical isomorphisms identify $\delta_1$ with
$(\pi\pi,\pi\rho,\rho\rho)$. Here the components are horizontal 
composites of the local product projection 2-cells. 
For example, $\pi\rho$ is
$$\bfig
\morphism(0,0)|a|/{@{->}@/^2em/}/<500,0>[A\ox A`A;p\wedge r]
\morphism(0,0)|b|/{@{->}@/_2em/}/<500,0>[A\ox A`A;r]
\morphism(500,0)|a|/{@{->}@/^2em/}/<500,0>[A`A\ox A;p^*\wedge r^*]
\morphism(500,0)|b|/{@{->}@/_2em/}/<500,0>[A`A\ox A;p^*]
\morphism(250,100)|r|<0,-200>[`;\rho]
\morphism(750,100)|l|<0,-200>[`;\pi]
\efig$$
We will write 
\begin{equation}\label{delta}
\delta=(\pi\pi,\pi\rho,\rho\rho)\f(p^*\wedge r^*)(p\wedge r)\ra%
p^*p\wedge p^*r\wedge r^*r\f A\ox A\ra A\ox A
\end{equation}
\eth
\prf
We have
$$p\wedge r\iso d^*(p\ox r)d\iso d^*(p,r)\iso d^*1_{A\ox A}=d^*$$
and 
$$p^*\wedge r^*\iso d^*(p^*\ox r^*)d\iso d^*(p\ox r)^*d\iso(p,r)^*d%
\iso1_{A\ox A}^*d=d$$
so that $dd^*\iso(p^*\wedge r^*)(p\wedge r)$. To exhibit 
the other isomorphism of the statement we will write
$d_3\f A\ox A\ra (A\ox A)\ox (A\ox A)\ox (A\ox A)$ for the 
three-fold diagonal map 
$(1_{A\ox A},1_{A\ox A},1_{A\ox A})$ and then
$$p^*p\wedge p^*r\wedge r^*r\iso d_3^*(p^*p\ox p^*r\ox r^*r)d_3%
\iso d_3^*(p^*\ox p^*\ox r^*)(p\ox r\ox r)d_3%
\iso(d^*\ox A)(A\ox d)$$
\frp

Of course $\delta=(\pi\pi,\pi\rho,\rho\rho)$ in (\ref{delta}) 
of the Lemma is invertible if and only if $A$ is Frobenius. 
We will write 
$$\nu=\rho\pi\f(p^*\wedge r^*)(p\wedge r)\ra r^*p\f A\ox A\ra A\ox A$$
for the ``other'' horizontal composite of projections and for $A$
Frobenius we define
$\mu$ as the unique 2-cell ($\nu.\delta\inv$) making commutative
\begin{equation}\label{mu}
\bfig
\Vtriangle[(p^*\wedge r^*)(p\wedge r)`p^*p\wedge p^*r\wedge r^*r`r^*p;%
\delta`\nu`\mu]
\efig
\end{equation}
We remark that a
local product of maps is not generally a map. (In the case of
the bicategory of relations a local product of maps is a partial
map.) Observe though that if $A$ is such that the maps 
$d\f A\ra A\ox A$ and $t\f A\ra I$ have right adjoints {\em in $\bM$}
then $A$ is a cartesian object in $\bM$ in the terminology of
\cite{ckw} and \cite{ckvw}. In this case $p\wedge r\f A\ox A\ra A$
is the map that provides ``internal'' binary products for $A$. 

For maps $f,g\f T{\s {\two}}A$ we write, as in \ref{pl3}, $A(f,g)$ 
for the composite $f^*g$ and
observe that the following three kinds of 2-cells are in natural 
bijective correspondence
$$\bfig
\morphism(0,0)|a|/{@{->}@/^2em/}/<500,0>[T`A;f]
\morphism(250,125)|m|<0,-250>[`;\alpha]
\morphism(0,0)|b|/{@{->}@/_2em/}/<500,0>[T`A;g]
\morphism(1000,0)|a|/{@{->}@/^2em/}/<500,0>[T`T;1_T]
\morphism(1250,125)|m|<0,-250>[`;\Hat\alpha]
\morphism(1000,0)|b|/{@{->}@/_2em/}/<500,0>[T`T;A(f,g)]
\morphism(2000,0)|a|/{@{->}@/^2em/}/<500,0>[T`A;g^*]
\morphism(2250,125)|m|<0,-250>[`;\alpha^*]
\morphism(2000,0)|b|/{@{->}@/_2em/}/<500,0>[T`A;f^*]
\efig$$
We have
\lem\label{homcat}
The hom-category $\bM(T,A)$ can be equivalently described as the
category whose objects are the maps $f\f T\ra A$ and whose hom-sets
$\bM(T,A)(f,g)$ are the sets $\bM(T,T)(1_T,A(f,g))$ with composition 
given by pasting composites of the form
$$\bfig
\morphism(1000,0)|a|/{@{->}@/^2em/}/<500,0>[T`T;1_T]
\morphism(1250,125)|m|<0,-250>[`;\Hat\alpha]
\morphism(1000,0)|b|/{@{->}@/_2em/}/<500,0>[T`T;A(f,g)]
\morphism(500,0)|a|/{@{->}@/^2em/}/<500,0>[T`T;1_T]
\morphism(750,125)|m|<0,-250>[`;\Hat\beta]
\morphism(500,0)|b|/{@{->}@/_2em/}/<500,0>[T`T;A(g,h)]
\morphism(500,0)|b|/{@{->}@/_6em/}/<1000,0>[T`T;A(f,h)]
\morphism(1000,-250)|m|<0,-250>[`;f^*\epsilon_gh]
\efig$$
\eth
\prf
It is a simple exercise with mates to show that the pasting composite
displayed is $\Hat{\beta\alpha}$. We note that $\Hat{1_f}=\eta_f$.
\frp

\lem\label{whisker}
For objects $f,h,g,k$ of $\bM(T,A)$, the whisker composite 
$$\bfig
\morphism(-500,0)|a|<500,0>[T`A\ox A;(g,k)]
\morphism(0,0)|a|/{@{->}@/^2em/}/<1000,0>%
[A\ox A`A\ox A;(p^*\wedge r^*)(p\wedge r)]
\morphism(0,0)|b|/{@{->}@/_2em/}/<1000,0>%
[A\ox A`A\ox A;p^*p\wedge p^*r\wedge r^*r]
\morphism(1000,0)|a|<500,0>[A\ox A`T;(f,h)^*]
\morphism(500,150)|m|<0,-300>[`;\delta=(\pi\pi,\pi\rho,\rho\rho)]
\efig$$
being in the notation of \ref{pl3}
$$(p^*\wedge r^*)(p\wedge r)((f,h)(g,k))\to^{\delta((f,h)(g,k))}%
(p^*p\wedge p^*r\wedge r^*r)((f,h)(g,k))$$
is
$$\bfig
\morphism(0,0)|a|<750,0>[T`A;g\wedge k]
\morphism(750,0)|a|<750,0>[A`T;f^*\wedge h^*]
\morphism(0,0)|b|/{@{->}@/_3em/}/<1500,0>%
[T`T;A(f,g)\wedge A(f,k)\wedge A(h,k)]
\morphism(750,-50)|m|<0,-200>[`;(\pi\pi,\pi\rho,\rho\rho)]
\efig$$
In fact, $(p\wedge r)(g,k)\iso g\wedge k$ and 
$(f,h)^*(p^*\wedge r^*)\iso f^*\wedge h^*$.
\eth
\prf
We have
$$(p\wedge r)(g,k)\iso d^*(p\ox r)d(g,k)\iso%
d^*(p\ox r)((g,k)\ox(g,k))d\iso d^*(g\ox k)d\iso g\wedge k$$
while
$$(f,h)^*(p^*\wedge r^*)\iso(f,h)^*d^*(p^*\ox r^*)d\iso%
((p\ox r)d(f,h))^*d\iso ((f\ox h)d)^*d\iso f^*\wedge h^*$$
On the other hand, precomposing with maps and postcomposing with
pams preserves local products so that we have
$$(f,h)^*(p^*p\wedge p^*r\wedge r^*r)(g,k)\iso%
(f,h)^*(p^*p)(g,k)\wedge(f,h)^*(p^*r)(g,k)\wedge(f,h)^*(r^*r)(g,k)$$
$$\iso f^*g\wedge f^*k \wedge h^*k$$
Assembling these results in hom-notation gives the statement.
\frp
The whisker composite in Lemma \ref{whisker} should be thought of
as the {\em instantiation} of $\delta$ at $((f,h)(g,k))$ and we
have been deliberately selective in mixing our notations in the
conluding diagram of the statement; $(\pi\pi,\pi\rho,\rho\rho)$
being more informative than $\delta((f,h)(g,k))$. If we
instantiate the rest of diagram (\ref{mu}) at $((f,h)(g,k))$,
which is to say whisker with $(f,h)^*(-)(g,k)$, then the result is
clearly the lower triangle below.
\begin{equation}\label{Xi}
\bfig
\Vtriangle(0,0)<750,500>%
[(f^*\wedge h^*)(g\wedge k)`A(f,g)\wedge A(f,k)\wedge A(h,k)`A(h,g);%
(\pi\pi,\pi\rho,\rho\rho)`\rho\pi`\mu((f,h),(g,k)]
\Atriangle(0,500)<750,500>%
[1_T`(f^*\wedge h^*)(g\wedge k)`A(f,g)\wedge A(f,k)\wedge A(h,k);%
\Xi`(\alpha,\beta,\gamma)`]
\efig
\end{equation}
In the top triangle above it is clear that a $1_T$-element of
$A(f,g)\wedge A(f,k)\wedge A(h,k)$ is exactly an ``S'' shaped 
configuration in $\bM(T,A)$ of the form
$$\bfig
\morphism(0,0)|a|<500,0>[f`g;\alpha]
\morphism(0,0)|m|<500,-500>[f`k;\beta]
\morphism(0,-500)|b|<500,0>[h`k;\gamma]
\efig$$
For $A$ Frobenius we will be interested in lifting $1_T$-elements of 
$A(f,g)\wedge A(f,k)\wedge A(h,k)$ though the isomorphism
$$(\pi\pi,\pi\rho,\rho\rho)\f(f^*\wedge h^*)(g\wedge k)\ra%
A(f,g)\wedge A(f,k)\wedge A(h,k)$$
As we discussed in \ref{pl3}, we do not have precise knowledge of
general $1_T$-elements $\Xi$ of 
$$(f^*\wedge h^*)(g\wedge k)=((p^*\wedge r^*)(p \wedge r))((f,h)(g,k))$$
but those obtained by pasting a $1_T$-element of
$(p^*\wedge r^*)((f,h),x)$ to a $1_T$-element of $(p\wedge r)(x,(g,k))$,
for some $x\f T\ra A$ present no difficulty. (Here, $p^*\wedge r^*$
is the $S$ and $p\wedge r$ is the $R$ of \ref{pl3}.) Since
$$(p^*\wedge r^*)((f,h),x)=(f,h)^*(p^*\wedge r^*)x\iso%
(f^*\wedge h^*)x\iso f^*x\wedge h^*x=A(f,x)\wedge A(h,x)$$
and
$$ (p\wedge r)(x,(g,k))=x^*(p\wedge r)(g,k)\iso%
x^*(g\wedge k)\iso x^*g\wedge x^*k=A(x,g)\wedge A(x,k)$$
(where we have used Lemma \ref{whisker} in each derivation) we see
that these special $1_T$-elements of $(f^*\wedge h^*)(g\wedge k)$
are given by (equivalence classes of)``X'' shaped configurations 
in $\bM(T,A)$ of the form
$$\bfig
\morphism(0,0)|a|<250,-250>[f`x;\xi]
\morphism(0,-500)|a|<250,250>[h`x;\zeta]
\morphism(250,-250)|b|<250,250>[x`g;\eta]
\morphism(250,-250)|b|<250,-250>[x`k;\omega]
\efig$$
It is convenient to write such a $1_T$-element of 
$(f^*\wedge h^*)(g\wedge k)$ as the following pasting composite
\begin{equation}\label{xconfig}
\bfig
\Atriangle(0,0)|lmb|[T`T`A;1_T`x`g\wedge k]
\Vtriangle(500,0)|amm|/->`->`<-/[T`T`A;1_T`x`x^*]
\Atriangle(1000,0)|mrb|/<-`->`->/[T`A`T;x^*`1_T`f^*\wedge h^*]
\morphism(500,350)|m|<0,-250>[`;(\eta,\omega)]
\morphism(1000,375)|m|<0,-250>[`;\eta_x]
\morphism(1500,350)|m|<0,-250>[`;(\xi^*,\zeta^*)]
\efig\end{equation}
Invertibility of 
$\delta=(\pi\pi,\pi\rho,\rho\rho)\f%
(f^*\wedge h^*)(g\wedge k)\ra X(f,g)\wedge X(f,k)\wedge X(h,k)$
tells us that, for every ``S'' configuration $(\alpha,\beta,\gamma)$,
there is a unique $1_T$-element $\Xi$ of $(f^*\wedge h^*)(g\wedge k)$
such that $\delta\Xi=(\alpha,\beta,\gamma)$. When, as in several
classical situations, every $1_T$-element $\Xi$ comes from an
``X'' configuration we have motivation for the colloquial name
``S''=``X'' for the Frobenius condition. (In fact one says
``S''=``X''=``Z'' when the second ``equation'' is not
derivable from the first but we have Lemma \ref{sym}.)

\lem\label{delX} For a $1_T$-element $\Xi$ (see (\ref{Xi})) arising 
from an ``X''configuration as in (\ref{xconfig}), 
$\delta\Xi=(\eta\xi,\omega\xi,\omega\zeta)$ and $\nu\Xi=\eta\zeta$.
\eth
\prf
For $\delta\Xi$ we treat the components separately. For the first, 
we paste 
$$\bfig
\morphism(0,0)|a|/{@{->}@/^2em/}/<500,0>[T`A;g\wedge k]
\morphism(0,0)|b|/{@{->}@/_2em/}/<500,0>[T`A;g]
\morphism(500,0)|a|/{@{->}@/^2em/}/<500,0>[A`T;f^*\wedge h^*]
\morphism(500,0)|b|/{@{->}@/_2em/}/<500,0>[A`T;f^*]
\morphism(250,100)|r|<0,-200>[`;\pi]
\morphism(750,100)|l|<0,-200>[`;\pi]
\efig$$
to (\ref{xconfig}) and obtain the $1_T$-element
$$\bfig
\Atriangle(0,0)|lmb|[T`T`A;1_T`x`g]
\Vtriangle(500,0)|amm|/->`->`<-/[T`T`A;1_T`x`x^*]
\Atriangle(1000,0)|mrb|/<-`->`->/[T`A`T;x^*`1_T`f^*]
\morphism(500,350)|m|<0,-250>[`;\eta]
\morphism(1000,375)|m|<0,-250>[`;\eta_x]
\morphism(1500,350)|m|<0,-250>[`;\xi^*]
\efig$$
of $A(f,g)$. To see this as a 2-cell $f\ra g$ paste onto it
$$\bfig
\morphism(0,0)|b|<500,0>[A`T;f^*]
\morphism(500,0)|a|<500,0>[T`A;f]
\morphism(0,0)|b|/{@{->}@/_3em/}/<1000,0>[A`A;1_A]
\morphism(500,-50)|r|<0,-200>[`;\epsilon_f]
\efig$$
(at $f^*$) which is the ``unhatting'' bijection and observe that the
result is $\eta\xi\f f\ra g$.
For the second, first paste $(\pi,\rho)$ and then paste $\epsilon_f$.
For the third, first paste $(\rho,\rho)$ and then paste $\epsilon_h$.
For $\nu\xi$, paste $(\rho,\pi)$ to (\ref{xconfig})
and then paste $\epsilon_h$ (at $h^*$).
\frp

The 2-cell $\mu$ of (\ref{mu}) when instantiated as in (\ref{Xi})
provides a completion of ``S'' configurations, as by the dotted arrow 
below. (It ultimately has the air of a Malcev operation.) 
$$\bfig
\morphism(0,0)|a|<500,0>[f`g;]
\morphism(0,0)|m|<500,-500>[f`k;]
\morphism(0,-500)|b|<500,0>[h`k;]
\morphism(0,-500)|b|/.>/<500,500>[h`g;]
\efig$$
In particular, given a 2-cell $\alpha\f f\ra g$ we have the
``S'' configuration $(1_f,\alpha,1_g)$ and we write
$\alpha^\dagger=\mu(1,\alpha,1)$.

\lem\label{main}
$\alpha^\dagger=\alpha\inv$
\eth
\prf
The composite $\alpha\alpha^\dagger$ is the clockwise composite
$1_T\ra A(g,g)$ in the following commutative diagram.
$$\bfig
\Ctriangle(0,0)|lml|/<-`->`->/%
[A(f,f)\wedge A(f,g)\wedge A(g,g)`1_T`A(f,g)\wedge A(f,g)\wedge A(g,g);%
(1,\alpha,1)`A(f,\alpha)\wedge A(f,g)\wedge A(g,g)`(\alpha,\alpha,1)]
\square(500,0)|amrb|<1000,1000>%
[A(f,f)\wedge A(f,g)\wedge A(g,g)`A(g,f)`%
A(f,g)\wedge A(f,g)\wedge A(g,g)`A(g,g);\mu`%
A(f,\alpha)\wedge A(f,g)\wedge A(g,g)`A(g,\alpha)`\mu]
\efig$$
We show that $\alpha\alpha^\dagger=1_g$ by evaluating the
counterclockwise composite. While we do not know if an ``X'' 
configuration gives rise to the $1_T$-element 
$\delta\inv(1,\alpha,1)$ we 
{\em do} know that $(\alpha,\alpha,1)$ arises from the
``X'' configuration
$$\bfig
\morphism(0,0)|a|<250,-250>[f`g;\alpha]
\morphism(0,-500)|a|<250,250>[g`g;1_g]
\morphism(250,-250)|b|<250,250>[g`g;1_g]
\morphism(250,-250)|b|<250,-250>[g`g;1_g]
\efig$$
because, writing $\Xi$ for the $1_T$-element arising as in 
(\ref{xconfig}) we have, by Lemma \ref{delX}, 
$\delta\Xi=(\alpha,\alpha,1)$. It follows using (\ref{Xi}) and again
Lemma \ref{delX}
that 
$$\alpha\alpha^\dagger=\mu(\alpha,\alpha,1)=\nu\Xi=1_g$$ 
Similarly, the composite $\alpha^\dagger\alpha$ 
is the clockwise composite in the commutative diagram.
$$\bfig
\Ctriangle(0,0)|lml|/<-`->`->/%
[A(f,f)\wedge A(f,g)\wedge A(g,g)`1_T`A(f,f)\wedge A(f,g)\wedge A(f,g);%
(1,\alpha,1)`A(f,f)\wedge A(f,g)\wedge A(\alpha,g)`(1,\alpha,\alpha)]
\square(500,0)|amrb|<1000,1000>%
[A(f,f)\wedge A(f,g)\wedge A(g,g)`A(g,f)`%
A(f,f)\wedge A(f,g)\wedge A(f,g)`A(f,f);\mu`%
A(f,f)\wedge A(f,g)\wedge A(\alpha,g)`A(\alpha,f)`\mu]
\efig$$
The rest of the proof proceeds as above after observing that
$(1,\alpha,\alpha)$ arises from the ``X'' configuration
$$\bfig
\morphism(0,0)|a|<250,-250>[f`f;1_f]
\morphism(0,-500)|a|<250,250>[f`f;1_f]
\morphism(250,-250)|b|<250,250>[f`f;1_f]
\morphism(250,-250)|b|<250,-250>[f`g;\alpha]
\efig$$
\frp

This completes the proof of Theorem \ref{grplhom}.


\end{document}